\documentclass [12pt,leqno]{article}

\newcounter{conjecture}
\newcounter{remark}
\newcounter{corollary}

\newenvironment{remark}{\medskip{\bf Remark \theremark.}
\addtocounter{remark}{1}}{}

\newcommand{\eqnsection}{
   \renewcommand{\theequation}{\thesection.\arabic{equation}}
   \makeatletter
   \csname @addtoreset\endcsname{equation}{section}
   \makeatother}
\newtheorem{theorem}{Theorem}
\newtheorem{lemma}{Lemma}

\def \be{\begin{equation}}
\def \ee{\end{equation}}
\def \bt{\begin{theorem}}
\def \et{\end{theorem}}
\def \bea{\begin{eqnarray}}
\def \eea{\end{eqnarray}}
\def \bas{\begin{eqnarray*}}
\def \eas{\end{eqnarray*}}

\def \al{\alpha}

\def \ga{\gamma}
\def \Ga{\Gamma} 
\def \de{\delta}
\def \De{\Delta}
\def \ep{\epsilon}

\def \la{\lambda}
\def \La{\Lambda}

\def \om{\omega}

\def \vf{\varphi}

\def \ze{\zeta}

\def \ff{\infty}

\def \rar{\rightarrow}

\def \R{{\bf R}}

\def \DD{{\cal D}}

\def \FF{{\cal F}}

\def \SS{{\cal S}}

\def \VV{{\cal V}}

\def \({\left(}
\def \){\right)}

\def \lc{\left\{}
\def \rc{\right\}}

\def \bsq{\ $\Box$}
\def \nn{\nonumber}

\def \bc{\begin{center} }
\def \ec{\end{center} }
\def \bs{\begin{slide} }
\def \es{\end{slide} }

\def\square{{\vcenter{\vbox{\hrule height.3pt
        \hbox{\vrule width.3pt height5pt \kern5pt
           \vrule width.3pt}
        \hrule height.3pt}}}}

\def\bsq{{\hfill $\square$ \bigskip}}

\eqnsection

\def \ilt{ intersection local time \,}

\def \VV{{\cal V}}
\def \nn{\nonumber}

\def \La{\Lambda}

\title{ Continuous Differentiability of Renormalized Intersection Local Times in $R^{1}$}
\author{Jay S. Rosen\thanks{This research   was   supported, in part, by grants from PSC-CUNY and the National Science
Foundation.}}
\date{}

\begin{document}
\maketitle

\renewcommand{\abstractname}{}
\begin{abstract}
 
Abstract.\hspace{.03 in}-\hspace{.03 in}We study $\gamma_{k}(x_2,\ldots,x_k;t)$, the k-fold
renormalized self-intersection local time
 for Brownian motion in $R^1$. Our main result says that  
$\gamma_{k}(x_2,\ldots,x_k;t)$  is   continuously differentiable in the spatial variables, with probability $1$.
 \vspace{.2in}

R\'esum\'e.\hspace{.03 in}-\hspace{.03 in}Nous \'etudions $\gamma_k(x_2,\ldots,x_k;t)$, le temps local  renormalis\'e d'auto-intersection  d'ordre k 
du mouvement brownien dans $ R^{1}$.
Notre r\'esultat principal montre que $\gamma_k(x_2,\ldots,x_k;t)$
est presque s\^urement contin\^ument diff\'erentiable dans les variables spatiales.
\end{abstract}


 
\eqnsection
\bibliographystyle{amsplain}

\section{Introduction}

The object of this paper is to establish the almost sure continuous differentiability of  renormalized intersection local time for the multiple intersections of Brownian motion in $R^{1}$.

Intersection local times were originally envisioned as a means of `measuring' the amount of
self-intersections of Brownian motion $W_{t} \in R^m$. Formally, the k-fold intersection local
time is
\[
\al_{k}(x_{2},x_{3},\ldots,x_{k}; t)=\int \cdots \int_{\{0\leq t_{1}\leq \cdots \leq t_{k}\leq t\}} \prod_{j=2}^{k} \delta
(W_{t_{j}}-W_{t_{j-1}}-x_{j})\,dt_{1}\ldots\,dt_{k}
\] where $\delta(x)$ denotes the $\delta$-function. Intuitively, $\al_{k}(0,0,\ldots,0; t)$ measures the  
`amount' of k-fold intersections.

More precisely,  we can set 
\bea
&&
 \al_{k,\ep}(x_{2},x_{3},\ldots,x_{k}; t) \label{0.1}\\
 &&=\int \cdots \int_{\{0\leq t_{1}\leq \cdots \leq t_{k}\leq t\}} \prod_{j=2}^{k}
f_{\ep}(W_{t_{j}}-W_{t_{j-1}}-x_{j})\,dt_{1}\ldots\,dt_{k}\nn
\eea where $f_{\ep}$ is an approximate $\delta$-function, and try to take the $\ep \rar 0$ limit.

In two dimensions, $\lim_{\ep\rar 0}\al_{k,\ep}(x_{2},x_{3},\ldots,x_{k}; t)$ will not exist unless all $x_{i}\neq 0$!  This gave rise to the problem of renormalization: to subtract from  $\al_{k,\ep}(x_{2},x_{3},\ldots,x_{k}; t)$ terms involving lower order intersection local times,
$\al_{j,\ep}(x_{2},x_{3},\ldots,x_{k}; t)$ for $j<k$, in such a way that a finite and continuous, $\ep\rar 0$ limit results. This was originally
done for double intersections of Brownian motion by Varadhan \cite{V}, and gave rise to
a large literature, see Bass and Khoshnevisan \cite{BK}, Dynkin \cite{D}, Le Gall \cite{LeGall4,LeGall-St.Flour,LeGall6}  and Rosen \cite{ JR1, JR2,JCDM,JC, ZQV}.  

In this paper we are concerned with one dimensional Brownian motion. In this case,  as we show below, the limit
\begin{equation}
\al_{k}(x_{2},x_{3},\ldots,x_{k}; t)=\lim_{\ep\rar 0} \al_{k,\ep}(x_{2},x_{3},\ldots,x_{k}; t).\label{0.2j}
\end{equation}
 exists a.s.  
Although, as  we will see, $\al_{k}(x_{2},x_{3},\ldots,x_{k}; t)$ is almost surely continuous, it is not 
$C^{1}$ in the spatial variable.  It is here that renormalization enters in the one dimensional case.

Let
\begin{equation}
 g(x)=\int_{0}^{\infty}e^{-t/2}p_{t}(x)\,dt=e^{-|x|}\label{gx}
\end{equation}
where $p_{t}(x)$ is the Brownian density function. 
 We
define the renormalized k-fold intersection local time 
 for $x=(x_{2},\ldots,x_{k})\in R ^{k-1}$ by
\be \ga_{k}(x;t)=\sum_{A\subseteq \{2,\ldots,k\}}(-1)^{|A|} (\prod_{j\in A}g(x_{j}) )
\al_{k-|A|}(x_{A^{c}};t) \label{0.11} \ee where for any
$B=\{i_{1}<\cdots<i_{|B|}\}\subseteq \{2,\ldots,k\}$ \be x_{B}=(x_{i_{1}}, x_{i_{2}}, \ldots,
x_{i_{|B|}})
 .\label{0.12} \ee 
 Here, we use the convention $\al_{1}(t)=t$. 
 Simple combinatorics then show that
\be
\al_{k}(x;t)=\sum_{A\subseteq \{2,\ldots,k\}} (\prod_{j\in A}g(x_{j}) )
\ga_{k-|A|}(x_{A^{c}};t).\label{com}
\ee
The renormalization (\ref{0.11}) used here is similar to that used in \cite{JC} and \cite{BK} for two dimensional Brownian motion, but in that case $g(x)=\int_{0}^{\infty}e^{-t/2}p_{t}(x)\,dt$ is infinite when $x=0$, compare (\ref{gx}). One key result of those papers is   that   $\ga_{k}(x;t)$ has a continuous extension from $(R^{2}-\{0\})^{k-1}\times R_{+}$ 
 to $(R^{2})^{k-1}\times R_{+}$.
   
Here is our main result.
 
\begin{theorem}\label{t1} For Brownian motion in $ R^1$ 
\be 
\al_{k}(x; t)=\lim_{\ep\rar 0}
\al_{k,\ep}(x; t) \label{5.1} 
\ee  exists and is jointly continuous a.s. Furthermore, $\ga_{k}(x; t)$ is  differentiable in $x$ and 
$\nabla_{x} \ga_{k}(x; t) $  is jointly continuous  with probability $1$.
\end{theorem}

For $k=2$, this was established in \cite{DSILT} by very different techniques. It seems impossible to use those techniques for $k>3$.

In \cite{JR3} we use  Theorem \ref{t1} to give a simple proof of the  CLT  for the $L^{2}$ modulus of continuity of  local time.

Note that  $ g(x)$ is continuously differentiable for $x\neq 0$. (\ref{com}) then exhibits precisely the non-differentiability of $\al_{k}(x;t)$. This justifies our choice of  renormalization (\ref{0.11}).  Simple combinatorics   show that we  obtain a similar result if we add any $C^{1}$ function to $ g(x)$. We have chosen $ g(x)$ as a potential density to simplify our proofs.

Our paper is organized along the lines of \cite{JC}. That paper was concerned with the continuity of $\ga_{k}(x_{2},x_{3},\ldots,x_{k}; t)$ in two dimensions, for Brownian motion and stable processes.   Our challenge here is to study differentiability, and for ease of exposition we consider only Brownian motion. After
laying the groundwork in Section $2$, we establish  the existence and almost sure continuity of $\al_{k}(x_{2},x_{3},\ldots,x_{k}; \ze)$
in Section
$3$, where $\ze$ is an independent mean-2 exponential random
variable. In Section
$4$ we show that the renormalized k-fold intersection local time $\ga_{k}(x_{2},x_{3},\ldots,x_{k}; \ze)$ is almost surely differentiable with a continuous derivative, again at an independent
exponential time. In Section
$5$ we use martingale techniques to obtain a.s. joint continuty.

\section{Intersection local times:  moments}

Let $W_{t}$ denote Brownian motion in $R^1$ with transition densites $p_{t}(x)$. In this
section we introduce approximate intersection local times for $W_t$ and present formulas for the
expectations of their moments. 

Let $f$ denote a smooth positive function supported on  $[-1,1]$,   such that  $\int
f(x)\,dx=1$, and for any
$\ep>0$ let
\[ f_{\ep}(y)=\frac{1}{\ep}\,f\(\frac{y}{\ep}\)
\] and $f_{\ep,x}(y)=f_{\ep}(y-x)$. We define the approximate \ilt of order $k$ as
\begin{eqnarray} & &
\al_{k,\ep}(x_{2},x_{3},\ldots,x_{k}; t) \label{1.1}\\ &=&\int \cdots \int_{\{0\leq t_{1}\leq
\cdots \leq t_{k}\leq t\}} \prod_{j=2}^{k} f_{\ep,x_{j}}(W_{t_{j}}-W_{t_{j-1}})\,dt_{1}\ldots\,dt_{k}. 
\nonumber
\end{eqnarray} 
We often abbreviate this as  $\al_{k,\ep}(x; t)$ where
$x=(x_{2},x_{3},\ldots,x_{k})\in R^{k-1}$. Let $g(x)=\int_{0}^{\infty}e^{-t/2}p_{t}(x)\,dt=e^{-|x|}$
denote the Green's function for $W_{t}$, and let $\ze$ denote a mean-$2$ exponential random
variable independent of $W_{t}$. The following   follows  appear in \cite[Section 2]{JC}, and are reproduced here for the convenience of the reader.  
The first formula follows easily from the Markov
property for $W_{t}$.

\begin{lemma}
\begin{eqnarray} & & E\(\prod_{i=1}^{n} \al_{k_{i},\ep_{i}}(x^{i}; \ze)\) \label{1.2}\\
&=&\sum_{v\in \VV} \int \prod_{i=1}^{n} \prod_{j=2}^{k_{i}}
 f_{\ep_{i},x_{j}^{i}}(y^{i}_{j}-y^{i}_{j-1}) \prod_{p=1}^{k}
 g(w_{v(p)}-w_{v(p-1)}) \,dw_{1} \ldots \,dw_{k}  \nonumber
\end{eqnarray} where $x^{i}=(x_{2}^{i},x_{3}^{i},\ldots,x_{k_{i}}^{i})$, $k=\sum_{i=1}^{n} k_{i}$,
$(w_{1},\ldots,w_{k})=(y^{1},\ldots,y^{n})\in (\R^{1})^{k}$ and $\VV$ is the set of permutations $v$ of
$\{1,2,\ldots,k\}$ such that whenever $ w_{v(p)}=y^{i}_{j}, w_{v(\tilde{p })}=y^{i}_{\tilde{j}}$ we
have $p>\tilde{p}\Longleftrightarrow j>\tilde{j}$.
\end{lemma}

A change of variables leads to the following more useful formula.

\begin{lemma}
\begin{eqnarray} & & E\(\prod_{i=1}^{n} \al_{k_{i},\ep_{i}}(x^{i}; \ze)\) \label{1.3}\\
&=&\sum_{s\in \SS} \int \prod_{i=1}^{n} \prod_{j=2}^{k_{i}}
 f(y^{i}_{j}) \prod_{p=1}^{k}
 g\(  z_{s(p)}+\sum_{j=2}^{c(p)}(\ep_{s(p)}y^{s(p)}_{j}+x^{s(p)}_{j})\right.
\nonumber \\ & &
-\left.(z_{s(p-1)}+\sum_{j=2}^{c(p-1)}(\ep_{s(p-1)}y^{s(p-1)}_{j}+x^{s(p-1)}_{j}))\)
\,dy^{i}_{j}\,dz_{1} \ldots \,dz_{n}  \nonumber
\end{eqnarray} where $x^{i}=(x_{2}^{i},x_{3}^{i},\ldots,x_{k_{i}}^{i})$, $k=\sum_{i=1}^{n} k_{i}$,
$\SS$ is the set of mappings $s:\{1,2,\ldots,k\}\mapsto \{1,\ldots,n\}$ such that
$|s^{-1}(i)|=k_{i},\forall 1\leq i\leq n$, and $c(p)=|\{u\leq p\mid s(u)=s(p)\}|$.
\end{lemma}

\section{Intersection local times: existence and continuity at exponential times}

We first consider the \ilt $\al_{k}$ at an independent mean-$2$ exponential time $\ze$.

A function $Z_{\ep}(x)$ indexed by $\ep\in (0,1]$ and $x$ in a topological space $\SS$ will be
said to converge locally uniformly in $x$ as $\ep\rar 0$ if for any compact $K\in \SS$, 
$Z_{\ep}(x)$ converges uniformly in $x\in K$ as $\ep\rar 0$.

\begin{theorem}\label{tilt} 
Almost surely, $ \al_{k,\ep}(x; \ze)$
converges locally uniformly in $x$ as $\ep\rar 0$.
Hence 
\be
\al_{k}(x; \ze):=\lim_{\ep\rar 0} \al_{k,\ep}(x; \ze) \label{3.2}
\ee 
 is continuous.

Furthermore, the occupation density formula holds:
\begin{eqnarray} & &
\int \Phi(x_{2},\ldots,x_{k}) \al_{k}(x_{2},\ldots,x_{k}; \ze)\,dx_{2}\ldots \,dx_{k}
\label{0.2b}\\ &=&\int
\cdots \int_{\{0\leq t_{1}\leq
\cdots \leq t_{k}\leq \ze\}} 
\Phi(W_{t_{2}}-W_{t_{1}},\ldots,W_{t_{k}}-W_{t_{k-1}})\,dt_{1}\ldots\,dt_{k} \nonumber
\end{eqnarray}
for all bounded Borel measurable functions $\Phi$ on $ R^{k-1}$.
\end{theorem}

\begin{remark} The occupation density formula (\ref{0.2b}) shows that $\al_{k}(x; \ze)$ is
independent of the particular $f$ used to define $\al_{k,\ep}(x; \ze)$.
\end{remark}

{\bf Proof of Theorem \ref{tilt}}: We will show that for $n$ even and $\ga>0$ we can find $\de>0$ such that 
\be E\( \{\al_{k,\ep}(x; \ze)-\al_{k,\ep'}(x'; \ze)\}^{n}\) \leq c_{n,\ga}|(\ep,x)-(\ep',x')|^{\de n/2} 
\label{3.3}
\ee for all $0<\ep,\ep' \leq \ga/2$ and all $x,x'\in R^{k-1}$. The multidimensional version of Kolmogorov's
Lemma, \cite[ Chapter 1, Theorem 2.1]{RY}, then gives us that for any  $\de'<\de$ and any $M<\infty$ we have
\be |\al_{k,\ep}(x; \ze)-\al_{k,\ep'}(x'; \ze)| \leq c_{n,\ga}(\om)|(\ep,x)-(\ep',x')|^{\de'/2} 
\label{3.4}
\ee for all rational $0<\ep,\ep' \leq \ga/2$ and all rational $x,x'\in R^{k-1}, |x|,|x'|
\leq M$. Since $\al_{k,\ep}(x; \ze)$ is clearly continuous as long as $\ep>0$, this will establish
the statements concerning (\ref{3.2}).

To establish (\ref{3.3}) we first handle the variation in $\ep$. If $h$ is a function of $\ep$, let \[\De_{\ep,\ep'}h=h(\ep)-h(\ep').\]
From Lemma 2 we have
\begin{eqnarray} & & E\(\prod_{i=1}^{n} \{ \al_{k ,\ep_{i}}(x^{i}; \ze)- \al_{k ,\ep'_{i}}(x^{i}; \ze)\}^{n}\) \label{3.5}\\
 & &= \prod_{i=1}^{n}\De_{\ep_{i},\ep'_{i}}   E\(\prod_{i=1}^{n} \al_{k ,\ep_{i}}(x^{i}; \ze)\)\nn \\
&&= \prod_{i=1}^{n}\De_{\ep_{i},\ep'_{i}} \sum_{s\in \SS} \int \prod_{i=1}^{n} \prod_{j=2}^{k }
 f(y^{i}_{j}) \prod_{p=1}^{nk}
 g\(  z_{s(p)}+\sum_{j=2}^{c(p)}(\ep_{s(p)}y^{s(p)}_{j}+x^{s(p)}_{j})\right.
\nonumber \\ & &
-\left.(z_{s(p-1)}+\sum_{j=2}^{c(p-1)}(\ep_{s(p-1)}y^{s(p-1)}_{j}+x^{s(p-1)}_{j}))\)
\,dy^{i}_{j}\,dz_{1} \ldots \,dz_{n}  \nonumber
\end{eqnarray} 
where we eventually set all $(\ep_{i},\ep'_{i})=(\ep,\ep')$. We expand this as a sum of many terms using
\be
 \De_{\ep,\ep'}(uv)= (\De_{\ep,\ep'}u)v(\ep)+ u(\ep')(\De_{\ep,\ep'}v)  \label{3.6}
\ee so that  each term  
contains for each $1\leq i\leq n$ a single difference of the form  $\De_{\ep_{i},\ep'_{i}} g$. Since each $g$ factor in
(\ref{3.5}) involves at most two $i$'s,
 whenever our procedure gives two differences involving the same $g$ factor we write one of
the differences as two terms. The upshot is that after setting all  $(\ep_{i},\ep'_{i})=(\ep,\ep')$, the expectation $E\( \{\al_{k,\ep}(x; \ze)-\al_{k,\ep'}(x; \ze)\}^{n}\)$ can be written as a
sum of many terms of the form appearing in  (\ref{1.3}) except that at least $n/2$ of the $g$
factors have been replaced by factors of the form
\bea
&&
\Delta_{\ep,\ep',j}
 g\Big (z_{s(p)}+\sum_{j=2}^{c(p)}(\tilde{\ep}_{s(p)}y^{s(p)}_{j}+x^{s(p)}_{j})\nonumber\\
 &&\hspace{1 in} -(z_{s(p-1)}+\sum_{j=2}^{c(p-1)}(\tilde{\ep}_{s(p-1)}y^{s(p-1)}_{j}+x^{s(p-1)}_{j}))\Big )
\label{3.7}
\eea  where $\tilde{\ep}$ can be variously $\ep,\ep'$ and the notation $\Delta_{\ep,\ep',j}  $ 
denotes a difference between two $g$ factors of the above form in which one of the $\ep$'s,
multiplying $y^{s(p)}_{j}$ or
$y^{s(p-1)}_{j}$ has been replaced by $\ep'$.

Our result then follows using 
\begin{equation}
|g(x)-g(y)|=|e^{-|x|}-e^{-|y|}|\leq |x-y|(g(x)+g(y))\label{lif.1}
\end{equation} for the variation in
$\ep$, and the variation in $x$ is handled similarly. 
This completes the proof of (\ref{3.2}).

To prove the occupation density formula (\ref{0.2b}) we note that
\begin{eqnarray} & &
\int \Phi(x_{2},\ldots,x_{k}) \al_{k,\ep}(x_{2},\ldots,x_{k}; \ze)\,dx_{2}\ldots \,dx_{k}
\label{0.2c}\\ &=&\int
\cdots \int_{\{0\leq t_{1}\leq
\cdots \leq t_{k}\leq \ze\}} 
\Phi\ast F_{\ep}(W_{t_{2}}-W_{t_{1}},\ldots,W_{t_{k}}-W_{t_{k-1}})\,dt_{1}\ldots\,dt_{k} \nonumber
\end{eqnarray}
where $F_{\ep}(x_2,\ldots,x_k)=\prod_{j=2}^k f_{\ep}(x_j)$. Hence, by what we have established
above, we can take the $\ep\rar 0$ limit in (\ref{0.2c}) to yield (\ref{0.2b}) whenever $\Phi$ is a
bounded continuous function. The monotone
convergence theorem then allows us to obtain (\ref{0.2b}) for all bounded Borel measurable
$\Phi$. This completes the proof of Theorem 2.\bsq

\section{Renormalized intersection local times:  continuous differentiability  at exponential times}

 We have 
defined the renormalized k-fold intersection local time 
 for   $x=(x_{2},\ldots,x_{k})\in R ^{k-1}$ by
\be \ga_{k}(x;t)=\sum_{A\subseteq \{2,\ldots,k\}}(-1)^{|A|} (\prod_{j\in A}g(x_{j}) )
\al_{k-|A|}(x_{A^{c}};t) \label{0.11cc} \ee where for any
$B=\{i_{1}<\cdots<i_{|B|}\}\subseteq \{2,\ldots,k\}$ \be x_{B}=(x_{i_{1}}, x_{i_{2}}, \ldots,
x_{i_{|B|}})
 .\label{0.12} \ee 
 Here, we use the convention $\al_{1}(t)=t$. 
Define the approximate renormalized k-fold intersection local time 
for $x=(x_{2},\ldots,x_{k})\in  R^{k-1}$ by
\be \ga_{k,\ep}(x;t)=\sum_{A\subseteq
\{2,\ldots,k\}}(-1)^{|A|} (\prod_{j\in A}g_{\ep}(x_{j}) )
\al_{k-|A|,\ep}(x_{A^{c}};t). \label{0.11acc} \ee
where
$g_{\ep}(x)=f_{\ep}\ast g(x)$. Clearly $\ga_{k,\ep}(x;\ze)\in C^{1}$ for fixed $\ep>0$ so that for any $y_{l},z_{l}$
\bea
&&
\ga_{k,\ep}(x_{2},\ldots,x_{l-1}, z_{l},x_{l+1},\ldots,x_{k } ;\ze)\label{4.0}\\
&&\hspace{.5 in}-\ga_{k,\ep}(x_{2},\ldots,x_{l-1}, y_{l},x_{l+1},\ldots,x_{k } ;\ze)=\int_{y_{l}}^{z_{l}} {\partial \over \partial x_{l}}\ga_{k,\ep}(x;\ze)\,dx_{l}.\nn
\eea

By Theorem \ref{tilt} and the continuity of $g$ it follows that almost surely 
$\ga_{k,\ep}(x;\ze)\rar \ga_{k}(x;\ze)$ as $\ep\rar 0$, locally uniformly in $x$. It follows from the next Theorem that ${\partial \over \partial x_{l}}\ga_{k,\ep}(x;\ze)\rar \psi_{k}(x;\ze)$ as $\ep\rar 0$, locally uniformly in $x$, for some continuous $\psi_{k}(x)$. Hence it follows from (\ref{4.0}) that
for any $y_{l},z_{l}$
\bea
&&
\ga_{k }(x_{2},\ldots,x_{l-1}, z_{l},x_{l+1},\ldots,x_{k} ;\ze)\label{4.0a}\\
&&\hspace{.5 in}-\ga_{k }(x_{2},\ldots,x_{l-1}, y_{l},x_{l+1},\ldots,x_{k} ;\ze)=\int_{y_{l}}^{z_{l}}  \psi_{k }(x;\ze)\,dx_{l}.\nn
\eea
This implies that $\ga_{k}(x;\ze)$ is differentiable with respect to $x_{l}$, and 
${\partial \over \partial x_{l}}\ga_{k}(x;\ze)= \psi_{k}(x;\ze)$, hence is continuous in $x$.  

\begin{theorem}\label{rilt} Almost surely, for each $2\leq l\leq k$
\be
{\partial \over \partial x_{l}}\ga_{k,\ep}(x;\ze)\label{3.2x}
\ee 
converges locally uniformly as $\ep\rar 0$ .
Hence the limit is continuous in $x$.
\end{theorem}

{\bf  Proof of Thorem \ref{rilt}: } As in the proof of Theorem \ref{tilt} it suffices to  show that for $n$ even and $\ga>0$ we can find $\de>0$ such that 
\be E\( \lc{\partial \over \partial x_{l}}\ga_{k,\ep}(x; \ze)-{\partial \over \partial x_{l}}\ga_{k,\ep'}(x'; \ze)\rc^{n}\) \leq c_{n,\ga}|(\ep,x)-(\ep',x')|^{\de n/2} 
\label{3.3k}
\ee for all $0<\ep,\ep' \leq \ga/2$ and all $x,x'\in R^{k-1}$.

Note that
\begin{eqnarray} & &
{\partial \over \partial x_{l}}\al_{k,\ep}(x_{2},x_{3},\ldots,x_{k}; t) \label{1.1der}\\ &=&-\int \cdots \int_{\{0\leq t_{1}\leq
\cdots \leq t_{k}\leq t\}} \prod_{j=2}^{l-1} f_{\ep}(W_{t_{j}}-W_{t_{j-1}}-x_{j})\nonumber\\
&& (f_{\ep})'(W_{t_{l}}-W_{t_{l-1}}-x_{l})\prod_{j=l+1}^{k} f_{\ep}(W_{t_{j}}-W_{t_{j-1}}-x_{j})\,dt_{1}\ldots\,dt_{k}. 
\nonumber
\end{eqnarray}
Note also  that since $g(x)$ is differentiable for all $x\neq 0$, with $g'(x)=-g(x)$ for $x>0$ and 
$g'(x)=g(x)$ for $x<0$, it follows that for any compactly supported $f\in C^{1}(R^{1})$
\begin{eqnarray}
&&
\int f'(x)g(x)\,dx
\label{ibp.1}\\
&& =\lim_{\ep\rar 0}\(\int_{-\ff}^{-\ep} f'(x)g(x)\,dx+\int_{-\ep}^{\ep} f'(x)g(x)\,dx+\int_{\ep}^{\ff} f'(x)g(x)\,dx\)
 \nonumber\\
&& =\lim_{\ep\rar 0}\(  f(-\ep)g(-\ep) - \int_{-\ff}^{-\ep} f(x)g(x)\,dx\)+\lim_{\ep\rar 0}\int_{-\ep}^{\ep} f'(x)g(x)\,dx
 \nonumber\\
&&\hspace{1 in} +\lim_{\ep\rar 0}\(  -f(\ep)g(\ep)+ \int_{\ep}^{\ff} f(x)g(x)\,dx\)\nonumber\\
&&= - \int_{-\ff}^{0} f(x)g(x)\,dx+ \int_{0}^{\ff} f(x)g(x)\,dx
 \nonumber\\
&&= - \int_{-\ff}^{\ff} f(x)g'(x)\,dx,  \nonumber
\end{eqnarray}
where for definiteness we set $g'(0)=0$.
Hence as in (\ref{1.3}), and using the product rule for differentiation
\begin{eqnarray} & & E\( \prod_{i=1}^{n} {\partial \over \partial x^{i}_{l}}\al_{k_{i},\ep_{i}}(x^{i}; \ze) \) \label{3.5j}\\
&&= (-1)^{n}  \sum_{s\in \SS} \int \prod_{i=1}^{n}\( \prod_{j=2,j\neq l}^{k_{i}}
 f_{\ep_{i}}(y^{i}_{j})\) \,(f_{\ep_{l}})'(y^{i}_{l})\nonumber\\
 &&
 \prod_{p=1}^{k}g\(  z_{s(p)}+\sum_{j=2}^{c(p)}(y^{s(p)}_{j}+x^{s(p)}_{j})\right.
\nonumber \\ & &
-\left.(z_{s(p-1)}+\sum_{j=2}^{c(p-1)}(y^{s(p-1)}_{j}+x^{s(p-1)}_{j}))\)
\,dy^{i}_{j}\,dz_{1} \ldots \,dz_{n}  \nonumber\\
&&=    \sum_{s\in \SS, a\in \mathcal{A}}(-1)^{\bar a_{2}} \int \prod_{i=1}^{n}  \prod_{j=2}^{k_{i}}
 f_{\ep_{i}}(y^{i}_{j})\nonumber\\
 &&  \hspace{.8 in}
 \prod_{p=1}^{nk}g^{(a_{1}(p)+ a_{2}(p))}\(  z_{s(p)}+\sum_{j=2}^{c(p)}(y^{s(p)}_{j}+x^{s(p)}_{j})\right.
\nonumber \\ & &\hspace{1 in}
-\left.(z_{s(p-1)}+\sum_{j=2}^{c(p-1)}y^{s(p-1)}_{j}+x^{s(p-1)}_{j}))\)
\,dy^{i}_{j}\,dz_{1} \ldots \,dz_{n}  \nonumber
\eea
where  $g^{(0)}=g,g^{(1)}=g',g^{(2)}=g''$, and $\mathcal{A}$ is the set of maps  $a=(a_{1}, a_{2}):\,[1,\ldots,kn]\mapsto \{0,1\}\times \{0,1\}$ such that 
\begin{itemize}
\item $\sum_{p=1}^{nk}a_{1}(p)+ a_{2}(p) =n$
\item for each 
$1\leq i\leq n$
\[\sum_{p=1}^{nk}a_{1}(p)1_{\{s(p)=i\}}+ a_{2}(p)1_{\{s(p-1)=i\}} =1 \]
\item if $a_{1}(p)=1$, then $c(p)\geq l$
\item  if $a_{2}(p)=1$, then $c(p-1)\geq l$.
\end{itemize}
In other words, if $s(p)=i$ then $a_{1}(p)=1$  if and only if,  after using the product rule for differentiation, the $p$'th $g$ factor is the only $g$ factor to which ${\partial \over \partial y^{i}_{l}}$
has been applied. Similarly, if $s(p-1)=i$ then $a_{2}(p)=1$  if and only if,  after using the product rule for differentiation, the $p$'th $g$ factor is the only $g$ factor to which ${\partial \over \partial y^{i}_{l}}$
has been applied. In (\ref{3.5j}), $\bar a_{2}=\sum_{p=1}^{nk}  a_{2}(p)$.

Then by scaling
\begin{eqnarray} & & E\( \prod_{i=1}^{n} {\partial \over \partial x^{i}_{l}}\al_{k_{i},\ep_{i}}(x^{i}; \ze) \) \label{3.5k}\\
&&=  \sum_{s\in \SS, a\in \mathcal{A}}(-1)^{\bar a_{2}}  \int \prod_{i=1}^{n} \prod_{j=2}^{k_{i}}
 f(y^{i}_{j})\nonumber\\
 &&\hspace{.3 in}
 \prod_{p=1}^{k}g^{(a_{1}(p)+ a_{2}(p))}\(  z_{s(p)}+\sum_{j=2}^{c(p)}(\ep_{s(p)}y^{s(p)}_{j}+x^{s(p)}_{j})\right.
\nonumber \\ & &\hspace{.5 in}
-\left.(z_{s(p-1)}+\sum_{j=2}^{c(p-1)}(\ep_{s(p-1)}y^{s(p-1)}_{j}+x^{s(p-1)}_{j}))\)
\,dy^{i}_{j}\,dz_{1} \ldots \,dz_{n}  \nonumber
\end{eqnarray} 
where we eventually set all $ \ep_{i} = \ep $.

We note in particular that if $p$ is a `bad integer', i.e. $s(p)=s(p-1)$, the $g^{(m)}$ factor in the above
product has the form
\be g^{(m)}(\ep_{s(p)}y^{s(p)}_{c(p)}+x^{s(p)}_{c(p)}) \label{4.2}
\ee
and in this case $m=0$ or $1$.

Let us now analyze the changes which occur in (\ref{3.5k}) when we replace the factor
${\partial \over \partial x^{r}_{l}}\al_{k_{r},\ep_{r}}(x^{r}; \ze)$ by ${\partial \over \partial x^{r}_{l}}\lc(\prod_{j\in B}g_{\ep_{r}}(x^{r}_{j}) ) \al_{k_{r}-|B|,\ep_{r}}(x^{r}_{B^{c}};\ze)\rc$.
Keeping in mind (\ref{4.2}) we see that now $s$ runs over those $s\in
\SS$ such that $s(p)=r, c(p)\in B\Rightarrow s(p-1)=r$, i.e. such $p$'s are bad, and in the
integrand on the right hand side of (\ref{3.5k}), aside from the factor $\prod_{j \in B}g^{(\cdot)}(\ep_{r}y^{r}_{j}+x^{r}_{j})$,
all other occurences of
$\ep_{r}y^{r}_{i}+x^{r}_{i},\,i\in B$ are deleted.

If $h(x)$ is any function of the variable $x$ we use the notation \[
\DD_{x}h =h(x)-h(0) \] for the difference between the value of $h$ at
$x$ and it's value at $x=0$. If $s\in \SS$ we set
$B_{s}=\{p|s(p)=s(p-1)\}$. The upshot is that we have 
\begin{eqnarray} & & E\( \prod_{i=1}^{n} {\partial \over \partial x^{i}_{l}}\ga_{k_{i},\ep}(x^{i}; \ze) \) \label{3.5m}\\
&&=   \sum_{s\in \SS, a\in \mathcal{A}}(-1)^{\bar a_{2}} \int \prod_{i=1}^{n} \prod_{j=2}^{k_{i}}
 f(y^{i}_{j})\(\prod_{p\in B_{s}}g^{(a_{1}(p)+ a_{2}(p))}(\ep_{s(p)}y^{s(p)}_{c(p)}+x^{s(p)}_{c(p)})\)\nonumber\\
 &&\(\prod_{p\in B_{s}}\DD_{\ep_{s(p)}y^{s(p)}_{j}+x^{s(p)}_{c(p)}}\)
 \prod_{p\in B_{s}^{c}}g^{(a_{1}(p)+ a_{2}(p))}\(  z_{s(p)}+\sum_{j=2}^{c(p)}(\ep_{s(p)}y^{s(p)}_{j}+x^{s(p)}_{j})\right.
\nonumber \\ & &\hspace{.5 in}
-\left.(z_{s(p-1)}+\sum_{j=2}^{c(p-1)}(\ep_{s(p-1)}y^{s(p-1)}_{j}+x^{s(p-1)}_{j}))\)
\,dy^{i}_{j}\,dz_{1} \ldots \,dz_{n}  \nonumber
\end{eqnarray} 
where we eventually set all $(\ep_{i},\ep'_{i})=(\ep,\ep')$.

To establish (\ref{3.3k}) we first handle the variation in $\ep$. By (\ref{3.5m})
\begin{eqnarray} & & E\( \prod_{i=1}^{n}\lc{\partial \over \partial x^{i}_{l}}\ga_{k_{i},\ep}(x^{i}; \ze)-{\partial \over \partial x^{i}_{l}}\ga_{k_{i},\ep'}(x^{i}; \ze)\rc \) \label{3.5mm}\\
&&= \prod_{i=1}^{n}\De_{\ep_{i},\ep'_{i}} \sum_{s\in \SS, a\in \mathcal{A}} (-1)^{\bar a_{2}}\int \prod_{i=1}^{n} \prod_{j=2}^{k_{i}}
 f(y^{i}_{j})\nonumber\\
 &&\(\prod_{p\in B_{s}}g^{(a_{1}(p)+ a_{2}(p))}(\ep_{s(p)}y^{s(p)}_{c(p)}+x^{s(p)}_{c(p)})\)
 \(\prod_{p\in B_{s}}\DD_{\ep_{s(p)}y^{s(p)}_{j}+x^{s(p)}_{c(p)}}\)\nonumber\\
 &&
 \prod_{p\in B_{s}^{c}}g^{(a_{1}(p)+ a_{2}(p))}\(  z_{s(p)}+\sum_{j=2}^{c(p)}(\ep_{s(p)}y^{s(p)}_{j}+x^{s(p)}_{j})\right.
\nonumber \\ & &
-\left.(z_{s(p-1)}+\sum_{j=2}^{c(p-1)}(\ep_{s(p-1)}y^{s(p-1)}_{j}+x^{s(p-1)}_{j}))\)
\,dy^{i}_{j}\,dz_{1} \ldots \,dz_{n}  \nonumber
\end{eqnarray} 
where we eventually set all $(\ep_{i},\ep'_{i})=(\ep,\ep')$. We will show that this is bounded in absolute value by
$c_{n}|\ep -\ep'|^{\de n/2}$  for all $|\ep|,|\ep'| \leq \ga$.  (At this stage, and for ease in generalizing to the variation in $x$, we allow 
  $\ep,\ep'$ to be zero or negative.) Compared to our proof of Theorem \ref{tilt}, the main difficulty here comes from the fact that $g'(x),\,g''(x)$, although uniformly bounded,  are not continuous at $x=0$.   It is here that the operators $\prod_{p\in B_{s}}\DD_{\ep_{s(p)}y^{s(p)}_{j}+x^{s(p)}_{c(p)}}$   will play a  critical role.

In the following we let $ g^{\sharp}=g,g'$ or $g''$. Fix $s\in \SS, a\in \mathcal{A}$. As in the proof of Theorem \ref{tilt}, the corresponding  term on  the left hand side of  (\ref{3.5mm}) can be written as a
sum of many terms of the form appearing in  (\ref{1.3}) except that at least $n/2$ of the $g^{\sharp}$
factors have been replaced by factors of the form
\be
\Delta_{\ep_{i} ,\ep'_{i} }
 g^{\sharp} (z_{i}+ \ep_{i}y_{j}^{i}  -z_{i'}+b_{i})
\label{3.7m}
\ee 
or
 \bea
 &&
\Delta_{\ep_{i} ,\ep'_{i} }\(
 g^{\sharp} (  \ep_{i}y_{j}^{i}  +b_{i}  )\DD_{\ep_{i}y_{j}^{i}  +b_{i} }\)= \(
 g^{\sharp}  (  \ep_{i}y_{j}^{i}  +b_{i}  )-g^{\sharp}  (  \ep'_{i} y_{j}^{i}  +b_{i} )\)\DD_{\ep_{i}y_{j}^{i}  +b_{i} }\nonumber\\
 &&\hspace{1 in} + 
 g^{\sharp}  (   \ep'_{i}y_{j}^{i}  +b_{i}  )\(\DD_{\ep_{i}y_{j}^{i}  +b_{i} }-\DD_{\ep'_{i}y_{j}^{i}  +b_{i} }\).
\label{3.7p}
\eea   
 Furthermore, we can assume that any $i$
  which appears in (\ref{3.7p})  differs from  from any $i,i'$ that appears in (\ref{3.7m}). Otherwise, we simply write one of them as a difference of two terms and consider each separately.
  
  Consider first  
  \begin{eqnarray}
  &&\Delta_{\ep_{i} ,\ep'_{i} }
 g^{\sharp} (z_{i}+ \ep_{i}y_{j}^{i}  -z_{i'}+b_{i})
  \label{lif.7}\\
  && =\Delta_{\ep_{i} ,\ep'_{i} }
 g^{\sharp} (z_{i}+ \ep_{i}y_{j}^{i}  -z_{i'}+b_{i})1_{\{|z_{i}+ \ep_{i}y_{j}^{i}  -z_{i'}+b_{i}|\geq 4|\ep_{i}-\ep'_{i}||y_{j}^{i}|\}}  \nonumber\\
  && +\Delta_{\ep_{i} ,\ep'_{i} }
 g^{\sharp} (z_{i}+ \ep_{i}y_{j}^{i}  -z_{i'}+b_{i})1_{\{|z_{i}+ \ep_{i}y_{j}^{i}  -z_{i'}+b_{i}|\leq 4|\ep_{i}-\ep'_{i}||y_{j}^{i}|\}}.  \nonumber
  \end{eqnarray}
  For the first term since we are away from the discontinuity of $g^{\sharp}$, we use (\ref{lif.1}) to obtain  a factor of $C|\ep_{i}-\ep'_{i}|$, since $|y_{j}^{i}|\leq 1$, while for the second term we use the fact that $g^{\sharp}$ is bounded and hence the $dz_{i}$
 integral contributes a factor of $C|\ep_{i}-\ep'_{i}|$. In more detail, on the set $\{|z_{i}+ \ep_{i}y_{j}^{i}  -z_{i'}+b_{i}|\leq 4|\ep_{i}-\ep'_{i}||y_{j}^{i}|\}$, up to a bounded error, we can replace every occurrence of $z_{i}$ in a $g^{\sharp}$ factor by $z_{i'}$, and in particular we simply bound 
$\Delta_{\ep_{i} ,\ep'_{i} }
 g^{\sharp} (z_{i}+ \ep_{i}y_{j}^{i}  -z_{i'}+b_{i})$ by $2$. This eliminates any occurrence of $z_{i}$ except in $1_{\{|z_{i}+ \ep_{i}y_{j}^{i}  -z_{i'}+b_{i}|\leq 4|\ep_{i}-\ep'_{i}||y_{j}^{i}|\}}$ which we can write as 
 $1_{\{B(z_{i'}- \ep_{i}y_{j}^{i}  -b_{i} ) ,4|\ep_{i}-\ep'_{i}||y_{j}^{i}|)\}}(z_{i})$. We then do the  the $dz_{i}$
 integral, to obtain a bound  of $C|\ep_{i}-\ep'_{i}|$.

A similar analysis holds for the last term in
(\ref{3.7p}) since $\DD_{\ep_{i}y_{j}^{i}  +b_{i} }-\DD_{\ep'_{i}y_{j}^{i}  +b_{i} }=\Delta_{\ep_{i}y_{j}^{i}  +b_{i} ,\ep'_{i}y_{j}^{i}  +b_{i} }$, and we have not yet `used' the $dz_{i}$
 integral. Finally, for the term $ \(
 g^{\sharp}  (  \ep_{i}y_{j}^{i}  +b_{i}  )-g^{\sharp}  (  \ep'_{i} y_{j}^{i}  +b_{i} )\)\DD_{\ep_{i}y_{j}^{i}  +b_{i} }$ in (\ref{3.7p}), if $|\ep_{i}y_{j}^{i}  +b_{i}  |\geq |\ep_{i}-\ep'_{i}|$ then $\ep_{i}y_{j}^{i}  +b_{i}$ and $\ep'_{i}y_{j}^{i}  +b_{i} $ have the same sign, so  we can use (\ref{lif.1}) to obtain  a factor of $C|\ep_{i}-\ep'_{i}|$, while if $|\ep_{i}y_{j}^{i}  +b_{i}  |\leq |\ep_{i}-\ep'_{i}|$ we can use $\DD_{\ep_{i}y_{j}^{i}  +b_{i}  }$ and the $dz_{i}$
 integral to obtain a factor of $C|\ep_{i}y_{j}^{i}  +b_{i}  |\leq C|\ep_{i}-\ep'_{i}|$. \bsq

\section{Renormalized intersection local times:  joint continuity of the spatial derivative}

Recall the approximate k'th order renormalized intersection local
time. 
\be \ga_{k,\ep}(x;t)=\sum_{A\subseteq \{2,\ldots,k\}}(-1)^{|A|} 
(\prod_{j\in A}g_{\ep}(x_{j}))
\al_{k-|A|,\ep}(x_{A^{c}};t). \label{8.7} \ee

\begin{theorem}\label{tjc} Almost surely,  $\ga_{\ep,k}(x;t)$ and  ${\partial \over \partial x_{l}}\ga_{k,\ep}(x;t)$   converge 
locally uniformly on $ R^{k-1} \times  R_{+}$ as
$\ep\rar 0$.
Hence 
\be
\ga_{k}(x;t)\stackrel{def}{=}\lim_{\ep\rar 0}\ga_{\ep,k}(x;t)\label{8.18}
\ee
is differentiable in $x$ and  $\nabla_{x}\ga_{\ep,k}(x;t)$    is continuous in $(x,t)\in R^{k-1} \times  R_{+}$.
\end{theorem}

\underline{Proof:}
Let $Y_t$ denote Brownian motion killed at  an independent mean-2
exponential time $\ze$. From now on, $\ga_{\ep,k}(x;t)$ will be defined for the
process $Y_t$ in place of $W_t$. Using Fubini's theorem, it suffices as before to show that $\ga_{\ep,k}(x;t)$ and ${\partial \over \partial x_{l}}\ga_{\ep,k}(x;t)$
converge  locally uniformly on $ R^{k-1} \times  [0,\ze)$ as
$\ep\rar 0$ with probability $1$. We will focus on ${\partial \over \partial x_{l}}\ga_{\ep,k}(x;t)$, and leave the easier case of $\ga_{\ep,k}(x;t)$ to the reader. 

If $\SS$ is a subset of Euclidean
space we will say that $\{Z_{\ep}(x);\,(\ep,x)\in (0,1]\times \SS\}$ converges rationally locally
uniformly on
$\SS$ as 
$\ep\rar 0$ if for any compact $K\in \SS$, 
$Z_{\ep}(x)$ converges uniformly in $x\in K$ as $\ep\rar 0$ when restricted to dyadic rational
$x,\ep$. We note that since ${\partial \over \partial x_{l}}\ga_{\ep,k}(x;t)$ for
$\ep>0$ is continuous in $\ep,x,t$, saying that ${\partial \over \partial x_{l}}\ga_{\ep,k}(x;t)$
converges 
 locally uniformly or converges rationally locally uniformly on $ R^{k-1} \times 
[0,\ze)$ as $\ep\rar 0$ are equivalent.

We know from Theorem \ref{rilt} that ${\partial \over \partial x_{l}}\ga_{\ep,k}(x;\ff)$ converges locally uniformly on 
$ R^{k-1}$ as $\ep\rar 0$ with probability $1$. Using martingale techniques we will see
that the right continuous martingale
\[
\Ga_{k,\ep,l}(x;t)\stackrel{def}{=}E\lc{\partial \over \partial x_{l}}\ga_{k,\ep}(x;\ff)\,\Bigg |\,\FF_t\rc
\]
converges rationally locally uniformly on $ R^{k-1} \times 
R_{+}$ as $\ep\rar 0$ with probability $1$. $\Ga_{k,\ep,l}(x;t)$ is not the same as
${\partial \over \partial x_{l}}\ga_{\ep,k}(x;t)$, but we will see that they differ by terms of `lower order', and we will be able
to complete our proof by induction. Given   the tools we have developed so far in this paper,
the proof  is conceptually fairly straightforward, but in order to treat the
`lower order' terms systematically we need to introduce some notation. This we now proceed to
do. 

We first define the approximate k'th order generalized intersection local
time
\begin{eqnarray} & &
\al_{k,\ep}(x_{2},x_{3},\ldots,x_{k};\phi; t) \label{8.2}\\ &=& \int_{\{0\leq t_{1}\leq
\cdots \leq t_{k}\leq t\}} \prod_{j=2}^{k}
f_{\ep,x_{j}}(Y_{t_{j}}-Y_{t_{j-1}})\phi(Y_{t_k})\,dt_{1}\ldots\,dt_{k}.  \nonumber
\end{eqnarray}
and set
\[
\al_{k,\ep}(x_{2},x_{3},\ldots,x_{k};\phi)=\al_{k,\ep}(x_{2},x_{3},\ldots,x_{k};\phi; \ff).
\]
$\al_{k,\ep}(x_{2},x_{3},\ldots,x_{k};\phi)$ is the approximate k'th order generalized 
total intersection local
time. For ease of notation in later formulas, we also set
\be
\al_{1,\ep}(\phi;t)= \int_{\{0\leq t_{1}\leq t\}}  \phi(Y_{t_1})\,dt_{1} \label{el1}
\ee
and
\be
\al_{0,\ep}(\phi;t)=\int \phi(z)\,dz\label{el2}
\ee
although $\al_{1,\ep}(\phi;t)$ is independent of $\ep$ and  $\al_{0,\ep}(\phi;t)$ is  independent of $\ep,t$.

Observe that
\begin{eqnarray}
&& 
E\{\al_{k,\ep}(x_{2},x_{3},\ldots,x_{k};\phi)\,| \,\FF_t\} \label{8.3} \\
&&
=\al_{k,\ep}(x_{2},x_{3},\ldots,x_{k};\phi; t)\nonumber\\ 
&&\qquad
+\sum_{i=0}^{k-1}E\( 
\int_{\{0\leq t_{1}\leq
\cdots\leq t_{i}\leq t \leq t_{i+1}\cdots\leq t_{k}\}} \right.\nonumber\\
&&\hspace{1in}
\left.\prod_{j=2}^{k}
f_{\ep,x_{j}}(Y_{t_{j}}-Y_{t_{j-1}})\phi(Y_{t_k})\,dt_{1}\ldots\,dt_{k}\,\Bigg|\,\FF_t\)  \nonumber
\end{eqnarray}
and for $i\geq 1$
\begin{eqnarray}
&&E\( 
\int_{\{0\leq t_{1}\leq
\cdots\leq t_{i}\leq t \leq t_{i+1}\cdots\leq t_{k}\}} \right.\label{done.1}\\
&&\hspace{1in}
\left.\prod_{j=2}^{k}
f_{\ep,x_{j}}(Y_{t_{j}}-Y_{t_{j-1}})\phi(Y_{t_k})\,dt_{1}\ldots\,dt_{k}\,\Bigg|\,\FF_t\)\nn
\\
&&=
\int_{\{0\leq t_{1}\cdots 
\leq t_{i}\leq t \leq t_{i+1}\cdots\leq t_{k}\}} \prod_{j=2}^{i}
f_{\ep,x_{j}}(Y_{t_{j}}-Y_{t_{j-1}})\nn\\
&& 
E\(f_{\ep,x_{i+1}}\((Y_{t_{i+1}}-Y_{t})+(Y_{t}-Y_{t_{i}})\)\prod_{j=i+2}^{k}
f_{\ep,x_{j}}(Y_{t_{j}}-Y_{t_{j-1}})\right.\nonumber\\
&&\hspace{.2 in}\left. \phi\((Y_{t_{k}}-Y_{t_{k-1}})+\cdots+(Y_{t_{i+1}}-Y_{t})+Y_{t}\)\,|\,\FF_t\)\,dt_{1}\ldots\,dt_{k} 
\nn\\
&&=
\int_{\{0\leq t_{1}\cdots 
\leq t_{i}\leq t \}} \prod_{j=2}^{i}
f_{\ep,x_{j}}(Y_{t_{j}}-Y_{t_{j-1}})f_{\ep,x_{i+1}}\(z_{i+1}+(Y_{t}-Y_{t_{i}})\)\nn\\
&&\hspace{.3 in} \prod_{j=i+2}^k f_{\ep,x_{j}}(z_j)
\phi(z_{i+1}+\cdots+z_k+Y_{t})\prod_{j=i+1}^k g(z_j)\,dz_{j} \,\,\,dt_{1}\ldots\,dt_{i} 
\nn\\
&&= \al_{i,\ep}(x_{2},x_{3},\ldots,x_{i};\la_{k-i,\ep}[\phi;x_{i+1},\ldots,x_k;Y_t];t)  \nonumber
\end{eqnarray}
where
\bea
&&
\la_{k-i,\ep}[\phi;x_{i+1},\ldots,x_k;u](v)\label{8.4}\\
&&\qquad
=\int f_{\ep,x_{i+1}}(z_{i+1}+u-v)\prod_{j=i+2}^k 
f_{\ep,x_j}(z_j)\nn\\
&&\hspace{1in}
\phi(z_{i+1}+\cdots+z_k+u)\prod_{j=i+1}^k g(z_j)\,dz_{i+1}\ldots\,dz_k\nn\\
&&\qquad
=\int g(z_{i+1}+v-u)\prod_{j=i+2}^k  g(z_j)\nn\\
&&\hspace{1in}
\phi(z_{i+1}+\cdots+z_k+v)\prod_{j=i+1}^k f_{\ep,x_j}(z_j)\,dz_{i+1}\ldots\,dz_k.\nn
\eea
Similarly, for $i=0$ we have 
\begin{eqnarray}
&&E\( 
\int_{\{0\leq   t \leq t_{1}\cdots\leq t_{k}\}} \right.\label{done.1z}\\
&&\hspace{1in}
\left.\prod_{j=2}^{k}
f_{\ep,x_{j}}(Y_{t_{j}}-Y_{t_{j-1}})\phi(Y_{t_k})\,dt_{1}\ldots\,dt_{k}\,\Bigg |\,\FF_t\)\nn
\\
&&=
\int_{\{0\leq   t \leq   t_{1}\cdots\leq t_{k}\}}  \nn\\
&& 
E\( \prod_{j=2}^{k}
f_{\ep,x_{j}}(Y_{t_{j}}-Y_{t_{j-1}})\right.\nonumber\\
&&\hspace{.2 in}\left. \phi\((Y_{t_{k}}-Y_{t_{k-1}})+\cdots+(Y_{t_{1}}-Y_{t})+Y_{t}\)\,|\,\FF_t\)\,dt_{1}\ldots\,dt_{k} 
\nn\\
&&=
\int   \prod_{j=2}^k f_{\ep,x_{j}}(z_j)
\phi(z_{1}+\cdots+z_k+Y_{t})\prod_{j=1}^k g(z_j)\,dz_{j}  
\nn\\
&&=
\int\(   \int g(v)
\phi(v+z_2+\cdots+z_k+Y_{t})\prod_{j=2}^k f_{\ep,x_{j}}(z_j) g(z_j)\,dz_{j} \)\,dv 
\nn\\
&&= \al_{0,\ep}( \la_{k,\ep}[\phi;x_{2},\ldots,x_k;Y_t];t)  \nonumber
\end{eqnarray}
where, recall our convention (\ref{el2}),
\bea
&&
\la_{k,\ep}[\phi;x_{2},\ldots,x_k;u](v)\label{8.4z}\\
&&\qquad
=\int g( v-u)
\phi(z_{2}+\cdots+z_k+v)\prod_{j=2}^k f_{\ep,x_{j}}(z_j) g(z_j)\,dz_{j}.\nn
\eea
By   abuse of notation we can introduce a fictitious $x_{1}$, and letting $f_{\ep,x_{1}}(z_1)$ denote $\de (z_{1})$, the $\de$-function, we can write
 \bea
&&
\la_{k,\ep}[\phi;x_{2},\ldots,x_k;u](v)\label{8.4y}\\
&& 
=\int g(z_1+ v-u)\prod_{j=2}^k   g(z_j)
\phi(z_1+z_{2}+\cdots+z_k+v)\prod_{j=1}^k f_{\ep,x_{j}}(z_j)  \,dz_{j}.\nn
\eea
Finally setting
\begin{equation}
  \la_{k,\ep}[\phi;x_{1},\ldots,x_k;Y_t] =:  \la_{k,\ep}[\phi;x_{2},\ldots,x_k;Y_t] \label{}
\end{equation}
 this now takes the same form as (\ref{8.4}) with $i=0$.
Then we can write
\begin{eqnarray}
&& 
E\{\al_{k,\ep}(x_{2},x_{3},\ldots,x_{k};\phi)\,|\,\FF_t\} \label{8.3z} \\
&&
=\al_{k,\ep}(x_{2},x_{3},\ldots,x_{k};\phi; t)\nonumber\\ 
&&\qquad
+\sum_{i=0}^{k-1}  \al_{i,\ep}(x_{2},x_{3},\ldots,x_{i};\la_{k-i,\ep}[\phi;x_{i+1},\ldots,x_k;Y_t];t) . \nonumber
\end{eqnarray}

Setting 
\bea
&&
\la_{k-i}[\phi;z_{i+1},\ldots,z_k;u](v)\label{8.5}\\
&&
= g(z_{i+1}+v-u)\prod_{j=i+2}^k
g(z_j)\phi(z_{i+1}+\cdots+z_k+v)\nn
\eea
we have
\bea
&&
\la_{k-i,\ep}[\phi;x_{i+1},\ldots,x_k;u](v)\label{8.6}\\
&&
=\int \la_{k-i}[\phi;z_{i+1},\ldots,z_k;u](v)\prod_{j=i+1}^k
f_{\ep,x_j}(z_j)\,dz_{i+1}\ldots\,dz_k.\nn
\eea

We next define the approximate k'th order generalized renormalized intersection local
time
\be \ga_{k,\ep}(x;\phi;t)=\sum_{A\subseteq \{2,\ldots,k\}}(-1)^{|A|} 
(\prod_{j\in A}g_{\ep}(x_{j}))
\al_{k-|A|,\ep}(x_{A^{c}};\phi;t) \label{8.7a} \ee
and set
\[
\ga_{k,\ep}(x_{2},x_{3},\ldots,x_{k};\phi)=\ga_{k,\ep}(x_{2},x_{3},\ldots,x_{k};\phi; \ff).
\]
$\ga_{k,\ep}(x_{2},x_{3},\ldots,x_{k};\phi)$ is the approximate k'th order generalized total
renormalized intersection local
time. As before, for ease of notation in later formulas, we also set
\[
\ga_{1,\ep}(\phi;t)= \int_{\{0\leq t_{1}\leq t\}}  \phi(Y_{t_1})\,dt_{1} 
\]
and
\[
\ga_{0,\ep}(\phi;t)=\int \phi(z)\,dz.
\]
Using (\ref{8.3z}) we find that
\begin{eqnarray} & &
E\{\ga_{k,\ep}(x_{2},x_{3},\ldots,x_{k};\phi)\,|\,\FF_t\} \label{8.8ja}\\
&&
=\ga_{k,\ep}(x_{2},x_{3},\ldots,x_{k};\phi; t)\nn\\ 
&& 
+\sum_{A\subseteq \{2,\ldots,k\}}(-1)^{|A|} 
(\prod_{j\in A}g_{\ep}(x_{j}))\nonumber\\
&&\hspace{.5 in}\sum_{i=0}^{k-|A|-1}
\al_{i,\ep}(x_{A^{c}(1)},x_{A^{c}(2)},\ldots,x_{A^{c}(i-1)};\nonumber\\
&&\hspace{1.5 in}\la_{k-|A|-i,\ep}[\phi;x_{A^{c}(i)},\ldots,x_{A^{c}(k-|A|)};Y_t];t)
\nonumber
\end{eqnarray}
where  
\[A^{c}=\{A^{c}(1)<A^{c}(2)<\cdots<A^{c}(k-|A|)\}.\]  
We   reorganize this by   writing $A$ as the disjoint union of 
\[A_{i}=\{j\in A\,|\, j<A^{c}(i)\} \mbox{  and  }B_{i}=\{j\in A\,|\, j>A^{c}(i)\}\]
so that
\begin{eqnarray}
&&(-1)^{|A|} 
(\prod_{j\in A}g_{\ep}(x_{j}))\nonumber\\
&&\hspace{.5 in}\al_{i,\ep}(x_{A^{c}(1)},x_{A^{c}(2)},\ldots,x_{A^{c}(i-1)};\nonumber\\
&&\hspace{1.5 in}\la_{k-|A|-i,\ep}[\phi;x_{A^{c}(i)},\ldots,x_{A^{c}(k-|A|)};Y_t];t)
\nn\\
&&=(-1)^{|A_{i}|} 
(\prod_{j\in A_{i}}g_{\ep}(x_{j}))  \nonumber\\
&&\hspace{.3 in} \al_{i,\ep}(x_{A^{c}(1)},x_{A^{c}(2)},\ldots,x_{A^{c}(i-1)};\nonumber\\
&&\hspace{.7  in}
(-1)^{|B_{i}|} 
(\prod_{j\in B_{i}}g_{\ep}(x_{j}))\la_{k-|A|-i,\ep}[\phi;x_{A^{c}(i)},\ldots,x_{A^{c}(k-|A|)};Y_t];t).\nonumber
\end{eqnarray}
It is now easy to see that if we fix $0\leq l\leq k-1$ and sum in (\ref{8.8ja}) over all $A\subseteq \{2,\ldots,k\}$
with $A^{c}(i)=l+1$ we will obtain 
\begin{equation}
\ga_{l,\ep}(x_{2},x_{3},\ldots,x_{l};\La_{k-l,\ep}[\phi;x_{l+1},\ldots,x_k;Y_t];t)\label{8.8jka}
\end{equation}
where
\bea
&&
\La_{k-l,\ep}[\phi;x_{l+1},\ldots,x_k;u](v)\label{8.11}\\
&&
=\int \La_{k-l}[\phi;z_{l+1},\ldots,z_k;u](v)\prod_{j=l+1}^k
f_{\ep,x_j}(z_j)\,dz_{i+1}\ldots\,dz_k\nn
\eea
with
\bea
&&
\La_{k-l}[\phi;z_{l+1},\ldots,z_k;u](v)\label{8.12}\label{8.11k}\\
&&=\sum_{B\subseteq \{l+2,\ldots,k\}}(-1)^{|B |} 
(\prod_{j\in B }g(z_{j}))\la_{k-|B|-l,\ep}[\phi;x_{\{l+1,\ldots,k\}-B};u] (v)\nonumber\\
&&= g(z_{l+1}+v-u)\sum_{B\subseteq \{l+2,\ldots,k\}}(-1)^{|B |} 
(\prod_{j\in B }g(z_{j}))\prod_{j\in \{l+2,\ldots,k\}-B}
g(z_j)\nonumber\\
&&\hspace{2.5 in}\phi(z_{l+1}+\sum_{i\in \{l+2,\ldots,k\}-B}z_i+v)\nn\\
&&= g(z_{l+1}+v-u)\prod_{j=l+2}^k
g(z_j)
\nonumber\\
&&\hspace{1  in}\sum_{B\subseteq \{l+2,\ldots,k\}}(-1)^{|B |} \phi(z_{l+1}+\sum_{i\in \{l+2,\ldots,k\}-B}z_i+v)\nn\\
&&
= g(z_{l+1}+v-u)\prod_{j=l+2}^k
g(z_j)\DD_{z_j}\phi(z_{l+1}+\cdots +z_{k}+v).\nn
\eea
The reader can check that this is consistent with our conventions when $l=0$.
Combining the above we obtain
\begin{eqnarray} & &
E\{\ga_{k,\ep}(x_{2},x_{3},\ldots,x_{k};\phi)\,|\,\FF_t\} \label{8.8a}\\
&&
=\ga_{k,\ep}(x_{2},x_{3},\ldots,x_{k};\phi; t)\nn\\ 
&&\qquad
+\sum_{i=0}^{k-1}
\ga_{i,\ep}(x_{2},x_{3},\ldots,x_{i};\La_{k-i,\ep}[\phi;x_{i+1},\ldots,x_k;Y_t];t).
\nonumber
\end{eqnarray}
We will say that $\La_{k-i}[\phi;z_{i+1},\ldots,z_k;u]$ is obtained from $\phi$ by adjunction
of $z_{i+1},\ldots,z_k;u$.

In order to make the sequel easier to follow, we make some explanatory comments. We will use (\ref{8.8a}) inductively  to show that  almost surely ${\partial \over \partial x_{l}}\ga_{k,\ep}(x_{2},x_{3},\ldots,x_{k}; t)$ converges locally uniformly in $x,t$ as $\ep\rar 0$. The convergence of the conditional expectation will follow easily from martingale inequalities and the techniques we have used to obtain convergence at exponential times. When looking at the last line  in (\ref{8.8a})  we encounter something new, the presence of $Y_t$. Rather than try to deal with this directly we prove that almost surely \[{\partial \over \partial x_{l}}\ga_{i,\ep}(x_{2},x_{3},\ldots,x_{i};\La_{k-i,\ep}[\phi;x_{i+1},\ldots,x_k;u];t)\] converges locally uniformly in $x,u,t$ as $\ep\rar 0$. Using the fact that $Y$ is almost surely locally bounded, this will show that almost surely \[{\partial \over \partial x_{l}}\ga_{i,\ep}(x_{2},x_{3},\ldots,x_{i};\La_{k-i,\ep}[\phi;x_{i+1},\ldots,x_k;Y_t];t)\] converges locally uniformly in $x,t$ as $\ep\rar 0$ which will allow us to complete the induction step. This explains the presence of $u$ in (\ref{8.11})-(\ref{8.12}). We will refer to such $u$ as a `new' parameter, while the $\{x_{2},\ldots, x_{k}\}$ will be referred to as `old' parameters.

We will say that a function $\vf_{y_{1},\ldots,y_n}(v)$ is {\it an admissable function
of $v$} with auxiliary parameters $y_{1},\ldots,y_n$ if it can be written in the
form
\be
\vf_{y_{1},\ldots,y_n}(v)=
\prod_{j\in B_0} g(y_j)\DD_{y_j}\prod_{i=1}^p g(v+\sum_{v\in B_i} \pm
y_v)\label{8.15}
\ee
where $B_i\subseteq \{1,\ldots,n\},\,\forall i=0,1,\ldots,p$,   
and $\{1,\ldots,n\}=\bigcup_{i=1}^p B_i$. Here $p$ is an arbitrary positive integer. If
$\vf_{y_{1},\ldots,y_n}(v)$ is of the above form we will say that $\vf_{y_{1},\ldots,y_n}(v)$ is of
weight $|B_0|+p$. Note that the weight of $\vf_{y_{1},\ldots,y_n}(v)$ is the number of $g$
factors in (\ref{8.15}). We will also consider the function
$\vf(v)\equiv 1$ to be an admissable function of $v$ (of weight $0$ and with no auxiliary
parameters). 

If $\vf_{y_{1},\ldots,y_n}(v)$ is an admissable function
of $v$ with auxiliary parameters $y_{1},\ldots,y_n$ we will use the notation
$\vf_{y_{1},\ldots,y_n;\ep}(v)$ to denote the function in which some of the auxiliary variables
have been smoothed. More precisely, we will say that $\vf_{y_{1},\ldots,y_n;\ep}(v)$ is a
{\it totally $\ep$-smoothed version of $\vf_{y_{1},\ldots,y_n}(v)$} if
\bea
\label{8.17}\\
\vf_{y_{1},\ldots,y_n;\ep}(v)=\int \vf_{y_{1}+\ep z_1,\ldots,y_n+\ep z_n}(v)\prod_{i\in
A}f(z_i)\,dz_i \prod_{i\in A^c} \,d\mu_0(z_i)\nn
\eea 
 for some subset
$A\subseteq \{1,\ldots,n\}$ such that (with the notation of (\ref{8.15})) $B_0\subseteq A$ and
$B_i\bigcap A\neq \emptyset$ for all 
$i=1,\ldots,p$. In other words,  we require that each $g$ factor in (\ref{8.15}) contain at least one element of
the set $y_j,\,j\in A$. Here $\mu_0$ is the Dirac measure which puts unit mass at the origin. It
would be more precise to refer to the function defined in (\ref{8.17}) as
$\vf_{y_{1},\ldots,y_n;\ep,A}(v)$, but in order to avoid further cluttering of the notation, and
because the actual nature of the set $A$ will be irrelevant for us, we shall simply drop it from
the notation. The reader will note in the sequel that it is precisely the `old' parameters which are integrated against an $f$.

The next lemma assembles some facts about adjunction which follow easily from the definitions.

\begin{lemma}\label{ladj}
Let $\vf_{y_{1},\ldots,y_n}(v)$ be an
admissable function of $v$ of weight $q$ and auxiliary parameters $y_{1},\ldots,y_n$, and 
let $\La_{k-i}[\vf_{y_{1},\ldots,y_n};x_{i+1},\ldots,x_k;u]$ denote the function in (\ref{8.12}) obtained from
$\vf_{y_{1},\ldots,y_n}(v)$ by adjunction of $x_{i+1},\ldots,x_k;u$. Then:
\begin{enumerate}
\item $\La_{k-i}[\vf_{y_{1},\ldots,y_n};x_{i+1},\ldots,x_k;u](z)$
is an admissable function of $z$ of weight
$q+k-i$ and auxiliary parameters $y_{1},\ldots,y_n,x_{i+1},\ldots,x_k,u$.
\item If $\vf_{y_{1},\ldots,y_n;\ep}(v)$ is a
totally $\ep$-smoothed version of $\vf_{y_{1},\ldots,y_n}(v)$, then thefunction 
$\La_{k-i,\ep}[\vf_{y;\ep} ;x_{i+1},\ldots,x_k;u]$ defined in (\ref{8.11}) is a totally
$\ep$-smoothed version of
$\La_{k-i}[\vf_{y};x_{i+1},\ldots,x_k;u]$.
\end{enumerate}
\end{lemma}

In the following, the notation $\nabla\ga_{i,\ep}(x;\vf_{y;\ep} )$ will denote the gradient with respect to $x$ and $y$. In fact, we are not interested in differentiating with respect to `new' parameters, but to avoid excessive notation we consider them also.

The next lemma generalizes Theorems \ref{rilt}.

\begin{lemma}\label{la}
Let $\vf_{y}(z)$ be an admissable function of $z$ of weight $k-i$ and 
auxiliary parameters $y=(y_1,\ldots,y_j)$ and let $\vf_{y;\ep} (z)$ be a
totally $\ep$-smoothed version of $\vf_{y}(z)$. Then  there exists $\de>0$ such that for each $n$ and $M<\ff$ we can find 
$c_{n,M}<\ff$ such that  
\be
 E\( \lc \sup_{F_M}\frac{
|\nabla\ga_{i,\ep}(x;\vf_{y;\ep} )-\nabla\ga_{i,\ep'}(x';\vf_{y';\ep'} )|}{|(\ep,x,y)-(\ep',x',y')|^{\de
}}\rc ^{n}\)\leq c_{n,M}\label{8.10}
\ee 
where $\sup_{F_M}$ is taken over all dyadic rational pairs $(\ep,x,y)\neq (\ep',x',y')$ such that 
$0<\ep,\ep'\leq 1$ and
$|x|,|x'|,|y|,|y'|\leq M$. 
\end{lemma}

 {\bf Proof of Lemma \ref{la}:} According to  \cite[ Chapter 1, Theorem 2.1]{RY}, it
suffices to show that
 there exists $\de>0$ such that for each $n$ and $M<\ff$ we can find  $c_{n,M}<\ff$
such that
\be\qquad
 E\( 
|\nabla\ga_{i,\ep}(x;\vf_{y;\ep} )-\nabla\ga_{i,\ep'}(x';\vf_{y';\ep'} )|
^{n}\)
\leq c_{n,M} |(\ep,x,y)-(\ep',x',y')|^{\de n}\label{8.16}
\ee 
for all $(\ep,x,y), (\ep',x',y')$ such that 
$0<\ep,\ep'\leq 1$ and
$|x|,|x'|,|y|,|y'|\leq M$.
(\ref{8.16}) follows as in the proof of Theorem  \ref{rilt}. \bsq

 {\bf Proof of Theorem \ref{tjc} (continued):}

We will show by induction on $i=0,1,\ldots,k$   that  $\nabla\ga_{i,\ep}(x;\vf_{y;\ep} ;t)$
converges 
locally uniformly in $(x,y,t)\in R^{j+i-1} \times 
[0,\ze)$ as $\ep\rar 0$ for all admissable functions $\vf_{y}(z)$ of $z$ of weight $k-i$ and 
auxiliary parameters $y=(y_1,\ldots,y_j)$.  

The case $i=k$ and $\vf_{y}(z)\equiv 1$ will prove our
theorem.

Consider first the case of $i=0$. We have to show that if $\vf_{y}(z)$
 is an admissable function of $z$ of weight $k$ and 
auxiliary parameters $y=(y_1,\ldots,y_j)$ and $\vf_{y;\ep} (z)$ is a
totally $\ep$-smoothed version of $\vf_{y}(z)$, then both 
\[
\ga_{0,\ep}(x;\vf_{y;\ep} ;t)\equiv \int \vf_{y;\ep} (z)\,dz
\]
and
\[
{\partial \over \partial y_{l}}\ga_{0,\ep}(x;\vf_{y;\ep} ;t)\equiv \int {\partial \over \partial y_{l}}\vf_{y;\ep} (z)\,dz
\]
converge locally uniformly in $y\in R^j$ as $\ep\rar 0$. This follows as in the  proof of Theorem  \ref{rilt}.

Assume now that for all $p<i$, and for all admissable functions $\Phi_{y}(z)$ of $z$ of weight $k-p$
and auxiliary parameters $y=(y_1,\ldots,y_{j'})$ we have
that $\nabla\ga_{p,\ep}(x;\Phi_{y;\ep} ;t)$ converges 
locally uniformly in $(x,y,t)\in R^{j'+p-1} \times 
[0,\ze)$ as $\ep\rar 0$ for any totally $\ep$-smoothed version $\Phi_{y;\ep} (z)$ of $\Phi_{y}(z)$.
Let us show that if
$\vf_{y}(z)$ is an admissable functions of $z$ of weight
$k-i$ and auxiliary parameters $y=(y_1,\ldots,y_j)$, and $\vf_{y;\ep} (z)$ is a
totally $\ep$-smoothed version of $\vf_{y}(z)$, then $\nabla\ga_{i,\ep}(x;\vf_{y;\ep} ;t)$
converges 
locally uniformly in $(x,y,t)\in R^{j+i-1} \times 
[0,\ze)$ as $\ep\rar 0$. 

With $F_M$ as in Lemma \ref{la},
let $F^m_M;\,\,m=1,2,\ldots$ be an exhaustion of $F_M$ by a
sequence of finite symmetric subsets. (A set $F$ of pairs $(a,b)$ is symmetric if $(a,b)\in F
\Rightarrow (b,a)\in F$). Let us define the right continuous martingale
\be
\Ga_{i,\ep}(x;\vf_{y;\ep} ;t)=E\(\nabla\ga_{i,\ep}(x;\vf_{y;\ep} )\,|\,\FF_t\).\label{8.22}
\ee
 By \cite[Chapter II, Theorem 1.7 ]{RY} applied to the right continuous submartingale
\be
A^m_t=\sup_{F^m_M}\frac{
|\Ga_{i,\ep}(x;\vf_{y;\ep} ;t)-\Ga_{i,\ep'}(x';\vf_{y';\ep'} ;t)|}
{|(\ep,x,y)-(\ep',x',y')|^{\de
}}
\label{8.23}
\ee
 we have that
\bea
&&
 E\( \lc \sup_{t}\sup_{F^m_M}\frac{
|\Ga_{i,\ep}(x;\vf_{y;\ep} ;t)-\Ga_{i,\ep'}(x';\vf_{y';\ep'} ;t)|}{|(\ep,x,y)-(\ep',x',y')|
^{\de
}}\rc ^{n}\)\label{8.24}\\
&&\qquad
 \leq E\( \lc \sup_{F^m_M}\frac{
|\nabla\ga_{i,\ep}(x;\vf_{y;\ep} )-\nabla\ga_{i,\ep'}(x';\vf_{y';\ep'} )|}{|(\ep,x,y)-(\ep',x',y')|
^{\de
}}\rc ^{n}\)\nn\\
&&\qquad
 \leq E\( \lc \sup_{F_M}\frac{
|\nabla\ga_{i,\ep}(x;\vf_{y;\ep})-\nabla\ga_{i,\ep'}(x';\vf_{y';\ep'} )|}{|(\ep,x,y)-(\ep',x',y')|^{\de
}}\rc ^{n}\)\nn\\
&&\qquad
\leq c_{n,M}\nn
\eea
where the last line used Lemma \ref{la}.
Hence
\be
 E\( \lc \sup_{t}\sup_{F_M}\frac{
|\Ga_{i,\ep}(x;\vf_{y;\ep} ;t)-\Ga_{i,\ep'}(x';\vf_{y';\ep'} ;t)|}
{|(\ep,x,y)-(\ep',x',y')|^{\de
}}\rc ^{n}\) \leq c_{n,M}. \label{8.25}
\ee
In particular this shows that
\be
\sup_{t}\sup_{F_{1,M}}|\Ga_{i,\ep}(x;\vf_{y;\ep} ;t)-\Ga_{i,\ep'}(x;\vf_{y';\ep'} ;t)|
\leq
C(\om)|\ep-\ep'|^{\de }\label{8.30}
\ee
where $F_{1,M}$ denotes the set of dyadic rational $(x,y)\in R^{j+i-1}$ with $|x|,|y|\leq M$.
Thus,
$\Ga_{i,\ep}(x;\vf_{y;\ep} ;t)$  converges rationally locally uniformly on $ R^{j+i-1}
\times 
R_{+}$ as $\ep\rar 0$ with probability $1$. 

It is easy to see that   $\nabla\ga_{i,\ep}(x;\vf_{y;\ep} ;t)$ is
continuous in
$\ep,x,y,t$ for $\ep>0$. Thus, as with $\nabla\ga_{\ep,k}(x;t)$, saying that
$\nabla\ga_{i,\ep}(x;\vf_{y;\ep};t)$ converges 
 locally uniformly or converges rationally locally uniformly  as $\ep\rar 0$ are equivalent.
By (\ref{8.8a}) 
\begin{eqnarray} & &
 \Ga_{i,\ep}(x;\vf_{y;\ep};t)\label{8.8b}\\
&&
=\nabla\ga_{i,\ep}(x_{2},x_{3},\ldots,x_{i};\vf_{y;\ep}; t)\nn\\ 
&&\qquad
+\sum_{p=0}^{i-1}
\nabla\ga_{p,\ep}(x_{2},x_{3},\ldots,x_{p};\La_{i-p,\ep}[\vf_{y;\ep};x_{p+1},\ldots,x_i;Y_t];t)
\nonumber
\end{eqnarray}
Hence to show that $\nabla\ga_{i,\ep}(x;\vf_{y;\ep};t)$ converges locally uniformly
on $ R^{j+i-1}
\times 
[0,\ze)$ as $\ep\rar 0$ with probability $1$ it suffices to show that for each $p<i$
\[ 
\nabla\ga_{p,\ep}(x_{2},x_{3},\ldots,x_{p};\La_{i-p,\ep}[\vf_{y;\ep};x_{p+1},\ldots,x_i;Y_t];t)
\]
converges locally uniformly
on $ R^{j+i-1}
\times 
[0,\ze)$ as $\ep\rar 0$ with probability $1$. However, by Lemma \ref{ladj}, 
$\La_{i-p,\ep}[\vf_{y;\ep};x_{p+1},\ldots,x_i;u]$ is a totally $\ep$-smoothed version of
$\La_{i-p}[\vf_{y};x_{p+1},\ldots,x_i;u]$, and the latter is an admissable function of weight
$k-p$ with auxiliary variables $y,x_{p+1},\ldots,x_i,u$. Therefore, by our induction
assumption,
\[ 
\nabla\ga_{p,\ep}(x_{2},x_{3},\ldots,x_{p};\La_{i-p,\ep}[\vf_{y;\ep};x_{p+1},\ldots,x_i;u];t)
\]
converges locally uniformly
in $(x,y,u,t)\in R^{j+i}
\times 
[0,\ze)$ as $\ep\rar 0$ with probability $1$. Since $Y_t$ is locally bounded on $[0,\ze)$, this 
completes proof of Theorem \ref{tjc}.
\bsq

\newcommand{\bysame}{\leavevmode\hbox to3em{\hrulefill}\,}

\bigskip
\noindent
\begin{tabular}{lll} 
      & Jay Rosen \\
      & Department of Mathematics \\
     &College of Staten Island, CUNY \\
     &Staten Island, NY 10314 \\ &jrosen30@optimum.net  
\end{tabular}

\end{document}